\documentclass{amsart}

\usepackage{graphicx}

\newtheorem{theorem}{Theorem}%[section]
\newtheorem{lemma}{Lemma}
\newtheorem{corollary}{Corollary}
\newtheorem{claim}{Claim}
\newtheorem{prop}{Proposition}

\theoremstyle{definition}

\theoremstyle{remark}

\numberwithin{equation}{section}

\begin{document}

\title{Asymptotic enumeration and limit laws of planar graphs}

\author{Omer Gim\'{e}nez}
\address{Departament de Matem\`{a}tica Aplicada II, Universitat Polit\`{e}cnica de
Catalunya, Jordi Girona 1--3, 08034 Barcelona, Spain}
\email{Omer.Gimenez@upc.edu}

\thanks{Research supported in part by Projects BFM2001-2340 and MTM2004-01728.}

\author{Marc Noy}
\address{Departament de Matem\`{a}tica Aplicada II, Universitat Polit\`{e}cnica de
Catalunya, Jordi Girona 1--3, 08034 Barcelona, Spain}
\email{Marc.Noy@upc.edu}


\date{\today}


\begin{abstract}
%Let $g_n$ be the number of  planar graphs on $n$ labelled
%vertices. We prove the  estimate $g_n \sim g\cdot n^{-7/2}
%\gamma^n n!$, where $g$ and $\gamma \approx 27.2268777685$ are
%computable constants. Next we find limit laws for several
%parameters in random planar graphs. We show that the number of
%edges is asymptotically normal with mean $\sim \kappa n$ for a
%computable constant $\kappa \approx 2.2132652385$, and that the
%number of edges is concentrated around the expected value. We
%show that the number of components in a random planar graph is
%distributed asymptotically as a shifted Poisson law $1+P(\nu)$,
%where $\nu \approx 0.0374393664$ is a computable constant. In
%particular, the probability of a random planar graph being
%connected tends to $e^{-\nu} \approx 0.9632528214$. We also prove
%and generalize a conjecture of McDiarmid, Steger and Welsh on the
%expected number of isolated vertices in random planar graphs.

We present a complete analytic solution to the problem of counting
planar graphs. We prove an estimate $g_n \sim g\cdot n^{-7/2}
\gamma^n n!$ for the number $g_n$ of labelled planar graphs on $n$
vertices, where $\gamma$ and $g$ are explicit computable
constants. We  show that the number of edges in random planar
graphs is asymptotically normal with linear mean and variance and,
as a consequence,  the number of edges is sharply concentrated
around its expected value. Moreover we prove an estimate
$g(q)\cdot n^{-4}\gamma(q)^n n!$ for the number of planar graphs
with $n$ vertices and $\lfloor qn \rfloor$ edges, where
$\gamma(q)$ is an analytic function of $q$.
We also show that the number of connected components in a random
planar graph is distributed asymptotically as a shifted Poisson
law $1+P(\nu)$, where $\nu$ is an explicit constant.
Additional Gaussian and Poisson limit laws for random planar
graphs are derived.
The proofs are based on singularity analysis of generating
functions and on perturbation of singularities.
\end{abstract}

\maketitle

\section{Introduction and statement of results}

In this paper we obtain a precise asymptotic estimate for the
number of labelled planar graphs on $n$ vertices, and we establish
limit laws for several parameters in random labelled planar
graphs. In particular, we  show that the number of edges in random
planar graphs is asymptotically normal, and that the number of
connected components in a random planar graph is distributed
asymptotically as a shifted Poisson law. Additional Gaussian and
Poisson limit laws for random planar graphs are derived.

 From now on, unless stated otherwise, all graphs are labelled.
Recall that a graph is planar if it admits an embedding in the
sphere. We remark that we consider planar graphs as combinatorial
objects, without referring to a particular topological embedding.

Let $g_n$ be the number of planar graphs on $n$ vertices. A
superadditivity argument~\cite{dominic} shows that the following
limit exists:
  $$
    \gamma = \lim_{n\to\infty} \left(g_n/n!\right)^{1/n}.
  $$
Until recently, the constant $\gamma$ was known only within
certain bounds, namely
 $$
    26.18 < \gamma < 30.06.
 $$
The lower bound results from the work of Bender, Gao and Wormald
\cite{bender}. They show that, if $b_n$ is the number of
2-connected  planar graphs, then
 $$\lim_{n\to\infty}
\left(b_n/n!\right)^{1/n} \approx 26.18.
$$
Hence $\gamma$ is at least this value.

The upper bound is based on the fact that an \emph{unlabelled}
planar graph on $n$ vertices can be encoded with at most $\alpha
n$ bits for some constant $\alpha$. If this is the case then $g_n
\le 2^{\alpha n} n!$, and so $\gamma \le 2^{\alpha}$. The first
such result was obtained by Tur\'{a}n~\cite{turan} with the value
$\alpha = 12$. This has been improved over the years and presently
the best result is $\alpha \approx 4.91$, obtained by Bonichon et
al.~\cite{upper2}. Since $2^{4.91} \approx 30.06$, the upper bound
follows.

Recently the present authors \cite{viena} were able to obtain,
using numerical methods, the approximation $\gamma \approx
27.2268$. In this paper we determine $\gamma$ exactly as an
analytic expression. Moreover, we find a precise asymptotic
estimate for the number of planar graphs.

\begin{theorem}\label{th:main}
Let $g_n$ be the number of planar graphs on $n$ vertices. Then
\begin{equation}\label{eq:est}
g_n \sim g \cdot n^{-7/2} \gamma^n n!,
\end{equation}
where
%%$g \approx 0.4970043999 \cdot 10^{-5}$
$g \approx 0.4260938569 \cdot 10^{-5}$
and $\gamma \approx 27.2268777685$.
\end{theorem}
The constants in the last statement are completely determined as
analytic expressions in terms of elementary functions. The proof
of~Thoerem~\ref{th:main}, together with the expressions given in
the appendix, contain all the necessary details for determinining
the constants. This also applies to all the remaining constants
that appear in the paper.

As we show later, for the number $c_n$ of \emph{connected}
 planar graphs on $n$ vertices, we have the estimate
 $$
 c_n \sim c \cdot n^{-7/2} \gamma^n n!,
 $$
where $\gamma$ is as before and
%%$c \approx 0.4787408907 \cdot 10^{-5}$.
$c \approx 0.4104361100 \cdot 10^{-5}$.

The proof of Theorem~\ref{th:main} is based on singularity
analysis of generating functions; see~\cite{FO,FS}. Let $g_n, c_n$
and $b_n$ be as before. As we show in the next section, there are
two equations linking the exponential generating functions
 $$
 B(x) = \sum b_n x^n/n!,
 \quad C(x) = \sum c_n x^n/n!, \quad G(x) =
\sum g_n x^n/n!.
$$
The dominant singularity of $B(x)$ was determined in
\cite{bender}; we are able to obtain the dominant singularities of
$C(x)$ and $G(x)$, which are both equal to~$\rho = \gamma^{-1}$.

In Section~\ref{pre} we review the preliminaries needed for the
proof. In Section~\ref{sec:b} we find an explicit expression for
the generating function $B(x,y)$ of 2-connected planar graphs
counted according to the number of vertices and edges. This is a
key technical result in the paper, which allows us to obtain a
full bivariate singular expansion of $B(x,y)$.
% in Lemma~\ref{sing-B}.
 The explicit expression obtained for the
function $\beta(x,y,z,w)$ in the statement of Lemma~\ref{le:key}
suggests that we are in fact integrating a rational function. This
is indeed the case as we explain later.

In Section~\ref{sec:est} we determine expansions of $C(x)$ and
$G(x)$ of square-root type at the dominant singularity $\rho$, and
then we apply ``transfer theorems''~\cite{FO,FS} to obtain
estimates for $c_n$ and $g_n$.

The singular expansions of $C(x)$ and $G(x)$ can be extended to
the corresponding bivariate generating functions $C(x,y)$ and
$G(x,y)$ near $y=1$. This allows us to prove in
Section~\ref{sec:laws}, using perturbation of
singularities~\cite{FS}, a normal limit law for the number of
edges in random planar graphs. To our knowledge, this problem was
first posed in \cite{dvw}.

Throughout this paper, we say that a sequence of random variables
$X_n$ with mean $\mu_n$ and variance $\sigma^2_n$ has a
\emph{normal limit law} if the normalized variables $X_n^* = (X_n
- \mu_n)/\sigma_n$ converge in law to the standard normal
distribution $\mathcal{N}(0,1)$; convergence in law means, as
usual, point-wise convergence of the corresponding distribution
functions.

\begin{theorem}\label{th:edges}
Let $X_n$ denote the number of edges in a random planar graph with
$n$ vertices. Then $X_n$ is asymptotically normal and the mean
$\mu_n$ and variance $\sigma_n^2$ satisfy
\begin{equation}\label{mean-var}
\mu_n \sim \kappa n, \qquad \sigma_n^2 \sim \lambda n,
\end{equation}
where $\kappa \approx2.2132652385$ and $\lambda \approx
0.4303471697$. The same is true, with the same constants, for
\emph{connected} random planar graphs.
\end{theorem}

As a consequence, since $\sigma_n = o(\mu_n)$, the number of edges
is concentrated around its expected value; that is, for every
$\epsilon >0$ we have
 $$
    \hbox{Prob} \{|X_n - \kappa n| > \epsilon n\} \to 0,
  \qquad \hbox{as $n \to \infty$}.
  $$
Previously it had been proved that $ \hbox{Prob} \{X_n <  \alpha
n\} \to 0$ and $ \hbox{Prob} \{X_n >  \beta n\} \to 0$,  as $n \to
\infty$, for some constants $\alpha$ and $\beta$. The best values
achieved so far were $\alpha \approx 1.85$ (shown in \cite{gm},
improving upon \cite{dvw}) and $\beta \approx 2.44$ (shown in
\cite{upper2}, improving upon \cite{opt}). Theorem~\ref{th:edges}
shows that in fact there is only one constant that matters,
namely~$\kappa$.

%deviations from $\kappa n$ are expected to be only of order $
%\mathcal{O}(\sqrt n)$, not of order~$\mathcal{O}(n)$

The previous theorem shows convergence in distribution to the
normal law. However, in this setting it is often the case that one
can also prove a \emph{local} limit law, that is convergence to
the \emph{density} function of the normal law.
We prove such a local limit law and we derive large deviation
estimates for the number of edges in random planar graphs.
In the next statement, as
later in the paper, $g_{n,q}$ and $c_{n,q}$ denote respectively
the number of planar graphs with $n$ vertices and $q$ edges;
$\rho(y)$, $G_5(y)$ and $C_5(y)$ are computable analytic functions
to be introduced later.

\begin{theorem}\label{th:local}
Let $\mu$ be a fixed ratio in the open interval $(1,3)$. Take
$u>0$ such that $-u \rho'(u)/\rho(u) = \mu$. Then, as $n$ goes to
$\infty$,
\begin{equation}\label{eq:local}
g_{n,\lfloor \mu n \rfloor} \sim n! \, G_5(u) {\rho(u)^{-n} u^{-\lfloor
\mu n \rfloor} \over \sqrt{2 \pi n}\, \Gamma(-5/2) \, \sigma n^{7/2}},
\end{equation}
where
 $$
 \sigma^2 =-u^2 \frac{\rho''(u)}{\rho(u)} - u\frac{\rho'(u)}{\rho(u)}
 +u^2\frac{\rho'(u)^2}{\rho(u)^2}.
 $$

The same is true for the number of \emph{connected} planar graphs
$c_{n,\lfloor \mu n \rfloor}$ if we replace $G_5(u)$ by $C_5(u)$ in
%Equation
(\ref{eq:local}).
\end{theorem}

\begin{figure}
\begin{center}
\includegraphics[width=\textwidth]{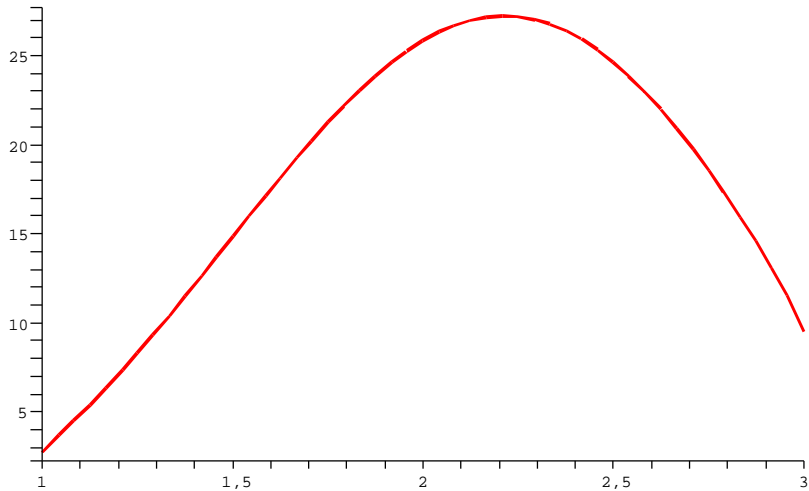}
\end{center}
\caption{The growth ratio of planar graphs with $n$ vertices and
$\lfloor \mu n \rfloor$ edges}\label{im:gamma}
\end{figure}

The previous result makes more precise a recent result from
\cite{gm2}, where the authors show that
 $$
    \lim_{n \to \infty} {1 \over n}  \log {g_{n,\lfloor \mu n \rfloor} \over n!}
     = \lambda(\mu),
 $$
where $\lambda(\mu)$ is a continuous function of $\mu$. A direct
consequence of Theorem \ref{th:local} is that
 $$
    \lambda(\mu)= -\mu \log(u)-\log(\rho(u))
 $$
where $u$ depends on $\mu$ as in the statement of the theorem.
Notice that $\lambda(\mu)$ is an analytic function of $\mu$.
Figure~\ref{im:gamma} shows the plot of $\exp(\lambda(\mu))$, that
is, the growth ratio of planar graphs with $n$ vertices and
$\lfloor \mu n \rfloor$ edges. The limit of $\exp(\lambda(\mu))$
as $\mu \to 1$ is equal to $e$, which is the growth ratio of
labelled trees; the limit as $\mu \to 3$ is equal to $256/27$,
which is the growth ratio of triangulations \cite{tutte-tri}.
(Tutte's result is for unlabelled triangulations, but a
triangulation has at most a linear number of automorphisms.)

Next we turn our attention to the following problem, considered in
\cite{dominic}.  Let $H$~be a graph on the vertex set
$\{1,\ldots,h\}$, and let $G$ be a graph on the vertex set
$\{1,\ldots,n\}$, where $n > h$.  Let $W \subset V (G)$ with $|W|
= h$, and let $r_W$ denote the least element in $W$. Following
\cite{dominic}, we say that $H$ \emph{appears} at $W$ in $G$ if
(a) the increasing bijection from $\{1,\ldots,h\}$ to $W$ gives an
isomorphism between $H$ and the induced subgraph $G[W]$ of $G$;
and (b) there is exactly one edge in $G$ between $W$ and the rest
of $G$, and this edge is incident with the root $r_W$.

Let $a_H (G)$ be the number of appearances of $H$ in $G$, that is,
the number of sets $W \subset V (G)$ such that $H$ appears at $W$
in $G$. Let $\alpha$ be $(9e^2(h+2))^{-1}\rho^h/h!$. It is shown
in \cite{dominic} that if $G_n$ is a random planar graph on $n$
vertices then
  $$
    {\rm Pr} \{a_H(G_n) \le \alpha n\} < e^{-\alpha n},
  $$
for $n$ large enough. The next result describes more precisely the
asymptotic behavior of the number of appearances of $H$ in random
planar graphs.

\begin{theorem}\label{th:appear}
Let $H$ be a fixed rooted connected planar graph with $h$
vertices. Let $X_n$ denote the number of appearances of $H$ in a
random planar graph with $n$ vertices. Then $X_n$ is
asymptotically normal and the mean $\mu_n$ and variance
    $\sigma_n^2$ satisfy
\begin{equation}\label{mean-appear}
\mu_n \sim {\rho^h \over h!}n, \qquad \sigma_n^2 \sim \rho n,
\end{equation}
where $\rho = \gamma^{-1}$ and $\gamma$ is as in Theorem
\ref{th:main}. Moreover, for every $\alpha<\rho^n/n!$ and every
$\beta>\rho^n/n!$ we have for $n$ large enough
\begin{equation}\label{prob-appear}
    {\rm Pr} \{X_n < \alpha n  \}<\left(
       \frac{u^{\alpha}}{x(u)\rho} \right)^n,
  \qquad
    {\rm Pr} \{X_n > \beta n  \}<\left(
       \frac{u^{\beta}}{x(u)\rho} \right)^n,
\end{equation}
where $x(u)$ is the solution of
$$
  xe^{(u-1)x^k/k!}=\rho,
$$
and $u$ is related to $Z$, where $Z$ is either $\alpha$ or
$\beta$, by the equation
$$
  -u\frac{x'(u)}{x(u)}=Z.
$$

\end{theorem}

Another parameter, the number of 2-connected components in a random
connected planar graph, also follows a normal distribution.

\begin{theorem}\label{th:blocs}
Let $X_n$ denote the number of blocks (2-connected components) in
a random \emph{connected} planar graph with $n$ vertices. Then
$X_n$ is asymptotically normal and the mean $\mu_n$ and variance
    $\sigma_n^2$ satisfy
\begin{equation}\label{mean-blocs}
\mu_n \sim \zeta n, \qquad \sigma_n^2 \sim \zeta n,
\end{equation}
where $\zeta \approx 0.0390518027$.
\end{theorem}

Next we turn to a different parameter, the number of connected
components in random planar graphs.

\begin{theorem}\label{th:components}
Let $X_n$ denote the number of connected components in a random
planar graph with $n$ vertices. Then  $X_n-1$ is distributed
asymptotically as a Poisson law of parameter $\nu$, where $\nu
\approx 0.0374393660$.
\end{theorem}

The above result is an improvement upon what was known so far. It
is shown in~\cite{dominic} that $Y_n$ is stochastically dominated
by $1+Y$, where $Y$ is a Poisson law $P(1)$;
Theorem~\ref{th:components} shows that in fact $Y_n$ is
asymptotically $1+P(\nu)$. The following direct corollary to
Theorem~\ref{th:components} is worth mentioning.

\begin{corollary}\label{th:coro}
\emph{(i)} The probability that a random planar graph is connected
is asymptotically equal to $e^{-\nu} \approx 0.9632528217$.
\emph{(ii)} The expected number of components in a random planar
graph is asymptotically equal to $1 + \nu \approx 1.0374393660$.
\end{corollary}

Our last result is the following. Let $\mathcal A$ be a family of
connected planar graphs, and let $A(x) = \sum A_n x^n /n!$ be the
corresponding generating function. Assume that the radius of
convergence of $A(x)$ is strictly larger than $\rho =
\gamma^{-1}$, the radius of convergence of $C(x)$; this is
equivalent to saying that $\mathcal A$ is exponentially smaller
than the family $\mathcal C$ of all connected planar graphs.

\begin{theorem}\label{th:family}
Assume $\mathcal A$ is a family of connected planar graphs that
satisfies the previous condition, and let $X_n$ denote the number
of connected components that belong to $\mathcal A$ in a random
planar graph with $n$ vertices. Then $X_n$ is distributed
asymptotically as a Poisson law of parameter $A(\rho)$.
\end{theorem}

If we take $\mathcal A$ as the family of graphs isomorphic to a
fixed connected planar graph~$H$ with $n$ vertices, then
 $$
 A(x)= {n! \over |\hbox{Aut}(H)|} \cdot {x^n \over n!} =
    {x^n \over |\hbox{Aut}(H)|},
 $$
where $\hbox{Aut}(H)$ is the group of automorphisms of $H$. In
particular, if $H$ is a single vertex, we obtain that the number
of isolated vertices in a random planar graph tends to a Poisson
law $P(\rho) = P(\gamma^{-1})$. This proves a conjecture by
McDiarmid, Steger and Welsh \cite{dominic}.

As a  different application of Theorem \ref{th:family} we have the
following. Recall that $B(x)$ is the generating function of
2-connected planar graphs.

\begin{corollary}
Let $X_n$ denote the number of connected components which are
2-connected in a random planar graph with $n$ vertices. Then $X_n$
tend to a Poisson law of parameter $B(\rho) \approx 0.0006837025$.
\end{corollary}

\medskip

We wish to emphasize that the approach that eventually has led to
the enumeration of planar graphs has a long history. Whitney's
theorem~\cite{whitney} guarantees that a 3-connected graph has a
unique embedding in the sphere; hence the problem of counting
3-connected graphs is in essence equivalent to counting
3-connected maps (planar graphs with a specific embedding). This
last problem was solved by Mullin and Schellenberg~\cite{MS} using
the approach developed by  Tutte in his seminal papers on counting
maps (see, for instance, \cite{tutte-maps}). The next piece is due
to Tutte~\cite{tutte}: a 2-connected graph decomposes uniquely
into 3-connected ``components''. Tutte's decomposition implies
equations connecting the generating functions of 3-connected and
2-connected planar graphs, which were obtained by
Walsh~\cite{walsh}, using the results of Trakhtenbrot \cite{trak}.
This in turn was used by Bender, Gao and Wormald~\cite{bender} to
solve the problem of counting 2-connected planar graphs; their
work is most relevant to us and is in fact the starting point of
our research. Finally, the decomposition of connected graphs into
2-connected components, and the decomposition of arbitrary graphs
into connected components, imply equations connecting the
corresponding generating functions. Analytic methods, together
with a certain amount of algebraic manipulation, become then the
main ingredients in our solution.

\emph{Acknowledgements.} We are grateful to Philippe Flajolet for
his encouragement and useful discussions during our research;  to
Eric Fusy for his help in deriving large deviation estimates and in
simplifying the final expressions in Lemma~\ref{le:key}; and to
Dominic Welsh for giving us access to an early version of
\cite{dominic}. Discussions with Manuel Bodirsky and Mihyun Kang
are also acknowledged.

%%%%%%%%%%%%%%%%%%%%%%%%%%%%%%%%%
\section{Preliminaries}\label{pre}
%%%%%%%%%%%%%%%%%%%%%%%%%%%%%%%%%

In this section and in the rest of the paper we use the language
and basic results of \emph{Analytic Combinatorics}, as in the
forthcoming book of Flajolet and Sedgewick~\cite{FS}.
% In particular, the theory of singular expansions and transfer
% theorems, and the extensions of the central limit theorem based on
% perturbation of singularities.
For the sake of completeness, we
state the main results we use in this paper (Corollary VI.1,
Theorems IX.10 and IX.13 in \cite{FS}).

\begin{prop}[Transfer Theorem; simplified version]
\label{prop:transfer} Assume that $f(z)$ is analytic in a domain
$\Delta=\Delta(\phi, R)$, where $R>1$, $0<\phi<\pi/2$ and
$$
  \Delta(\phi, R)=\{z\,:\, z\neq 1, |z|<R, \,|Arg(z-1)|>\phi\}.
$$
If, as $z\rightarrow 1$ in $\Delta$,
 $$
  f(z) \sim (1-z)^{-\alpha}
  $$
then
 $$
    [z^n]f(z) \sim \frac{n^{\alpha-1}}{\Gamma(\alpha)}.
 $$
\end{prop}

\begin{prop}[Quasi-Powers Theorem; algebraic singularities]
\label{th:quasi-powers}
Let $f(z,u)$ be a bivariate function that is bivariate analytic at $(0,0)$
with nonnegative coefficients there. Assume that $f$ admits in
$\mathcal{D}=\{|z|\leq r\}\times \{|u-1|<\epsilon\}$, for some $r>0$ and
$\epsilon>0$, the representation
$$
  f(z,u)=A(z,u)+B(z,u)C(z,u)^{-\alpha},
$$
where $A$, $B$ and $C$ are analytic in $\mathcal{D}$ such that
$C(z,1)=0$ has a unique simple root $\rho<r$ in $|z|\leq r$ and
$B(\rho,1)\neq 0$. Moreover, neither $\partial_z C(\rho,1)$ nor
$\partial_u C(\rho,1)$ are $0$, so there exists a nonconstant $\rho(u)$
analytic at $u=1$ such that $C(\rho(u),u)=0$ and $\rho=\rho(1)$.
Finally, $\rho(u)$ is such that
\begin{equation}\label{eq:quasi-powers}
  -\frac{\rho''(1)}{\rho(1)}-\frac{\rho'(1)}{\rho(1)}+\left(
\frac{\rho'(1)}{\rho(1)}\right)^2
\end{equation}
is non-zero. Then the random variable with probability generating
function
$$
  p_n(u)=\frac{[z^n]f(z,u)}{[z^n]f(z,1)}
$$
converges in distribution to a Gaussian variable. The mean $\mu_n$
and the standard deviation $\sigma_n$ converge asymptotically to
$\mu n$ and $\sigma \sqrt{n}$, where $\mu$ is $-\rho'(1)/\rho(1)$ and
$\sigma^2$ is given by (\ref{eq:quasi-powers}).
\end{prop}

\begin{prop}[Local Limit Law; simplified version]
\label{prop:local} Let $f(x,u)$ satisfy the conditions of
Proposition \ref{th:quasi-powers}.  If $\rho(u)$ attains uniquely
its minimum on the circle $|u|=1$ at $u=1$, then the sequence
$p_{n,\lfloor \mu n\rfloor}$ is asymptotically $(\sqrt{2 \pi
n}\,\sigma)^{-1}$.
\end{prop}

Now we discuss the generating functions that appear in this paper.
Recall that $g_n$, $c_n$ and $b_n$ denote, respectively, the
number of planar graphs, connected planar graphs, and 2-connected
planar graphs on $n$ vertices. The corresponding exponential
generating functions are related as follows.

\begin{lemma}\label{le:BCG}
The series $G(x)$, $C(x)$ and $B(x)$ satisfy the following
equations:
 $$
    G(x) = \exp(C(x)), \qquad xC'(x) =
    x\exp\left(B'(xC'(x))\right),
 $$
where $C'(x)$ and $B'(x)$ are derivatives with respect to $x$.
\end{lemma}

\begin{proof}
The first equation is standard, given the fact that a planar graph
is a \emph{set} of connected  planar graphs, and the set
construction in labelled structures corresponds to taking the
exponential of the corresponding exponential generating function.

The second equation follows from a standard argument on the
decomposition of a connected graph into 2-connected components.
Take a connected graph rooted at a vertex $v$; hence the
generating function $xC'(x)$. Now $v$ belongs to a \emph{set} of
2-connected components (including single edges), each of them
rooted at vertex $v$; hence the term $\exp(B')$. Finally, in each
of the 2-connected components, replace every vertex by a rooted
connected graph; this explains the substitution $B'(xC'(x))$.
Details can be found, for instance, in \cite[p. 10]{HP}.
\end{proof}

 Let $b_{n,q}$ be the number of 2-connected planar graphs with
$n$ vertices and $q$ edges, and let
 $$
    B(x,y) = \sum b_{n,q} y^q {x^n \over n!}
 $$
be the corresponding bivariate generating function. Notice that
$B(x,1) = B(x)$. The generating functions $C(x,y)$ and $G(x,y)$
are defined analogously. Since the parameter ``number of edges''
is additive under taking connected and 2-connected components, the
previous lemma can be extended as follows.

\begin{lemma}\label{le:BCG-biv}
The series $G(x,y)$, $C(x,y)$ and $B(x,y)$ satisfy the following
equations:
 $$
    G(x,y) = \exp(C(x,y)), \qquad x {\partial \over \partial x } C(x,y)    =
    x\exp\left({\partial \over \partial x} B(x {\partial \over \partial x } C(x,y),y)
    \right).
 $$
\end{lemma}

In the remaining of the section we recall the necessary
results from~\cite{bender}. Define the series $M(x,y)$ by means of
the expression
\begin{equation}\label{eq:M}
    M(x,y) = x^2y^2 \left({1 \over 1+xy} + {1 \over 1+y}-1-
    {(1+U)^2(1+V)^2 \over (1+U+V)^3} \right),
\end{equation}
where $U(x,y)$ and $V(x,y)$ are algebraic functions given by
\begin{equation}\label{eq:UV}
    U = xy(1+V)^2, \quad V = y(1+U)^2.
\end{equation}

In the next result and in the rest of the paper, all logarithms
are natural.

\begin{lemma}[Bender et al. \cite{bender}]\label{th:B}
We have
\begin{equation}\label{eq:B-D}
 {\partial B(x,y) \over \partial y} = {x^2 \over 2}
 \left(
 {1+D(x,y)  \over 1+y}
 \right),
\end{equation}
where $D = D(x,y)$ is defined implicitly by $D(x,0)=0$ and
\begin{equation}\label{eq:D}
  {M(x,D) \over 2x^2D} - \log\left({1+D \over 1+y}\right) +
  {xD^2 \over 1+xD} = 0.
\end{equation}
Moreover, the coefficients of $D(x,y)$ are nonnegative.
\end{lemma}
There is a small modification in equation (\ref{eq:B-D}) with
respect to~\cite{bender}. We must consider the graph consisting of
a single edge as being 2-connected, otherwise Lemmas~\ref{le:BCG}
and~\ref{le:BCG-biv} would not hold. Hence the term of lowest
degree in the series $B(x,y)$ is $y x^2/2$.

Let us comment on the previous equations. The algebraic generating
function~$M$ corresponds to (rooted) 3-connected planar maps.
%, andwas found by Mullin and Schellenberg \cite{MS}.
 The decomposition
of a 2-connected graph into 3-connected components implies
equations (\ref{eq:B-D}) and (\ref{eq:D}),
% which were found by Walsh \cite{walsh}.
The generating function $D(x,y)$ is that of
planar \emph{networks}, which are special graphs with two distinguished
vertices.

We define the following functions of the complex variable~$t$. The
appendix contains additional functions that are introduced later.

\begin{eqnarray*}
    \xi & = &  {(1+3t)(1-t)^3  \over 16t^3} \\
    Y & = & {(1+2t) \over (1+3t)(1-t)}\exp\left(-
    {t^2(1-t)(18+36t+5t^2) \over 2(3+t)(1+2t)(1+3t)^2}  \right) -1    \\
    \alpha & = & 144+592t+664t^2+135t^3+6t^4-5t^5 \\
    \beta & = & 3t(1+t)(400+1808t+2527t^2+1155t^3+237t^4+17t^5) \\
    D_0 & = & {3t^2 \over (1-t)(1+3t)} \\
    D_2 & = & - {48t^2(1+t)(1+2t)^2(18+6t+t^2) \over (1+3t)\beta} \\
    D_3 &=& 384t^3(1+t)^2(1+2t)^2(3+t)^2\alpha^{3/2}\beta^{-5/2}
\end{eqnarray*}
Let us notice a slight change in terminology: functions $\xi$ and
$Y$ are denoted, respectively, $x_0$ and $y_0$ in \cite{bender};
also, we correct a typo, namely a $t$ factor that was missing in
the expression for $D_2$.

A key fact is that for $y$ in a suitable small neighborhood of~1,
the equation $Y(t) = y$ has a unique solution in $t = t(y)$. Then
define
\begin{equation}\label{eq:R}
R(y) = \xi(t(y)). \end{equation}
 In the next lemma, $D_i(y)$
stands for $D_i(t(y))$. This applies too to functions $B_i(y)$ and
$C_i(y)$ that we introduce later in the paper.

\begin{lemma}[Bender et al. \cite{bender}]\label{le:sing}
For fixed $y$ in a small neighborhood of\/ $1$, $R(y)$ is the
unique dominant singularity of $D(x,y)$. Moreover, $D(x,y)$ has a
branch-point at $R(y)$, and the singular expansion at $R(y)$ is of
the form
$$
     D(x,y) = D_0(y) + D_2(y)X^2 + D_3(y)X^3 +
     \mathcal{O}(X^4),
 $$
where $X = \sqrt {1-x/R(y)}$ and the $D_i(y)$ are as before.
\end{lemma}

The previous lemma is the key result used in \cite{bender} to
prove the estimate
 $$
     b_n \sim b \cdot n^{-7/2} R^{-n} n!,
 $$
where $b$ is a constant and $R=R(1) \approx 0.0381910976$.

%%%%%%%%%%%%%%%%%%%%%%%%%%%%%%%%%
\section{Analysis of $B(x,y)$}\label{sec:b}
%%%%%%%%%%%%%%%%%%%%%%%%%%%%%%%%%

 From equation (\ref{eq:B-D}), it follows that
\begin{equation}\label{eq:int}
     B(x,y) = {x^2\over 2} \int_0^y
     {1+D(x,t) \over 1+t} \,dt.
\end{equation}
Our goal is to obtain an expression for $B(x,y)$ as a function of
$x,y$ and $D(x,y)$ that, although more complex, does \emph{not}
contain an integral. Recall that the algebraic function $U$ is
defined in (\ref{eq:UV}), and $D$ is defined in Lemma~\ref{th:B}.

\begin{lemma}\label{le:key}
Let $W(x,z) = z(1+U(x,z))$. The generating function $B(x,y)$ of
2-connected planar graphs admits the following expression as a
formal power series:
\begin{equation}\label{Bmain}
 B(x,y) =
\beta\left(x,y,D(x,y),W\left(x,D(x,y)\right)\right),
\end{equation}
where
$$
   \beta(x,y,z,w)=\frac{x^2}{2} \beta_1(x,y,z) -
      \frac{x}{4} \beta_2(x,z,w),
$$
and
\begin{eqnarray*}
 \beta_1(x,y,z) &=& {z(6x-2+xz) \over 4x} +
 (1+z)\log\left({1+y \over 1+z}\right) - {\log(1+z) \over 2}
 + {\log(1+xz) \over 2x^2}; \\
 \beta_2(x,z,w) &=& {2(1+x)(1+w)(z+w^2) + 3(w-z)  \over 2(1+w)^2}
  - {1 \over 2x} \log(1 + xz + xw + xw^2) \\
 &+&
  {1-4x \over 2x} \log(1+w) +
  {1-4x+2x^2 \over 4x} \log\left({1 - x +xz + -xw + xw^2 \over
 (1-x)(z+w^2+1+w)} \right). \\
\end{eqnarray*}

\end{lemma}

\begin{proof}
 From equation (\ref{eq:int}) we obtain
$$
     B(x,y) = {x^2\over 2}\log(1+y) +
    {x^2 \over 2} \int_0^y {D(x,t) \over 1+t} dt.
$$
We integrate by parts and obtain
 $$
   \int_0^y {D(x,t) \over 1+t} \, dt = \log(1+y) D(x,y) -
   \int_0^y \log(1+t) {\partial D(x,t) \over \partial t} \, dt.
  $$
 From now on $x$ is a fixed value. Now notice that from
(\ref{eq:D}) it follows that
 $$
    \phi(u) = -1 + (1+u)\exp\left( - {M(x,u) \over 2x^2u} - {xu^2
    \over 1+xu}\right),
 $$
is an inverse of $D(x,y)$, in the sense that $\phi(D(x,y)) = y$.
In the last integral we change variables $s= D(x,t)$, so that
$t=\phi(s)$. Then
\begin{eqnarray*}
    \int_0^y \log(1+t) {\partial D(x,t) \over \partial t} \, dt
    &=&    \int_0^{D(x,y)} \left(\log(1+s) - {xs^2 \over 1+xs}
    \right)ds  \\
    &-& \int_0^{D(x,y)} {M(x,s) \over 2x^2s} \, ds.
\end{eqnarray*}
The first integral has a simple primitive and we are left with an
integral involving $M(x,y)$.
Summing up we have
\begin{equation}\label{newB}
 B(x,y)= \Theta(x,y,D(x,y)) + {1 \over 4} \int_0^{D(x,y)}
 {M(x,s) \over s} \, ds,
\end{equation}
where $\Theta$ is the elementary function
 $$
 \Theta(x,y,z) = {x^2 \over 2} \left(z + {1 \over 2} z^2 %-y
 +   (1+z) \log{1+y \over 1+z} \right) -{x \over 2}z + {1
 \over 2} \log(1+xz).
 $$

Now we concentrate on the last integral. From (\ref{eq:M}) and
(\ref{eq:UV}) it follows that
 $$
    \int_0^D  {M(x,s) \over s} \, ds = x \int_0^D {(1+U)^2 U \over
    (1+U+V)^3} \, ds,
$$
where $U$ and $V$ are considered as functions of $x$ and $s$, and
where for simplicity we write $D = D(x,y)$ from now on.

 From the definition $W(x,s) = s(1+U(x,s))$, we obtain that
 $$
    {(1+U)^2 U \over (1+U+V)^3} = {W-s \over W(1+W)^3}.
    $$
Since $W$ satisfies the equation
 $$
        xs^2 + (1+2xW^2)s + W(xW^3-1) = 0,
$$
the functional inverse of $W(x,s)$ with respect to the second
variable is equal to
\begin{equation}\label{inverse}
    -t^2 - {1 - \sqrt{1+4xt+4xt^2} \over 2x},
\end{equation}
where we use  $t$ to denote the new variable.

It follows that
 $$
\int_0^D\frac{W-s}{W(1+W)3} \, ds = \int_0^{W(x,D)}
 \frac{\left(Q-1-2xt-2xt^2\right)\left(2Qt-2t-1\right)}
 {2xt(1+t)^3Q} \, dt,
 $$
where for simplicity we write
\begin{equation}\label{defQ}
 Q(x,t)=\sqrt{1+4xt+4xt^2}.
\end{equation}
The last integral can be solved explicitly with the help of a
computer algebra system such as {\sc Maple}, and we obtain as a
primitive the function

\begin{eqnarray*}
%&&
  \frac{1-2(t+4x+4xt)}{4x(1+t)^2}
 -\frac{1+2x(1+t)}{2x(1+t)^2}Q^3
 +\left(2+4xt+\frac{1+2(t-x-tx)}{4x(1+t)^2}\right)Q + \qquad
  \\
 %&&
 \frac{2x^2-4x+1}{4x} \log\left(
 \frac{Q+(1-2x-2xt)}{Q-(1-2x-2xt)}\right)
 -\frac{1}{2x}\log(Q+1+2xt)+\frac{1-4x}{2x}\log(1+t).
\end{eqnarray*}

Finally we have to replace $t$ for $W(x,D)$ in the previous
equation. The expression (\ref{inverse}) and equation (\ref{defQ})
imply that
 $$
 Q(x,W(x,D))=1+2x(D+W(x,D)^2).
$$
Hence when replacing $t$ for $W(x,D)$ we obtain an expression in
$x$, $D$ and $W(x,D)$ that is free of square roots. A routine
computation, combined with the intermediate equation (\ref{newB}),
gives the final expression for $B(x,y)$ as claimed.
\end{proof}

The function $\beta$ in the previous lemma looks like the
primitive of a rational function. This can be explained as follows
(we are grateful to P. Flajolet for this observation). The
algebraic equation satisfied by $U$ (here $x$ is considered as a
parameter)~is
 $$
 u - xy(1+(y(1+u)^2)^2 =0.
 $$
It can be checked (for instance, using the Maple package {\tt
algcurves}), that this equation in $u$ and $y$ defines a rational
curve, that is a curve of genus zero, and so it admits a rational
parametrization $(u(t),y(t))$. Now an integral $\int R(s,U(x,s))\,
ds$, where $R$ is a rational function, becomes the integral of a
rational function after the change of variables $s = y(t)$. In
particular, this applies to the integral in equation~(\ref{newB}).

The former lemma can be used to obtain the singular expansion of
$B(x,y)$. The function $R(y)$ is defined in (\ref{eq:R}) and $B_0,
B_2, B_4, B_5$ are analytic functions of $y$ given in the
appendix.  Again $B_i(y)$ stands for $B_i(t)$, where $t$ is the
unique solution of $Y(t)=y$ in a  neighborhood of~$1$.

\begin{lemma}\label{sing-B}
For fixed $y$ in a small neighborhood of\/ $1$, the dominant
singularity of $B(x,y)$ is equal to $R(y)$. The singular expansion
at $R(y)$ is of the form
 $$
 B(x,y) = B_0(y) + B_2(y) X^2 + B_4(y) X^4 + B_5(y) X^5 +
 \mathcal{O} (X^6),
 $$
where $X = \sqrt {1-x/R(y)}$, and the $B_i$ are analytic functions
in a neighborhood of~1.
\end{lemma}

\begin{proof}
Consider the  expression for $B(x,y)$ in Lemma~\ref{le:key} as a
function of $x,y$ and $D(x,y)$.
%We know from Pringsheim's theorem
%that a function with non-negative coefficients has a singularity
%at its radius of convergence, hence we focus
A simple analysis shows that, for $y$ close to $1$, the only
singularities come from the singularities of $D(x,y)$, hence the
first claim follows.

For the second assertion, substitute the singular expansions of
$D(x,y)$ and $U(x,D(x,y))$ (taken, respectively, from
Lemma~\ref{le:sing} and the appendix)  for $D(x,y)$ and
$U(x,D(x,y))$ in~(\ref{Bmain}) (recall that $W=z(1+U)$). Next set
$x = \xi(t)(1-X^2)$ and $y=Y(t)$ as functions of $t$, and expand
the resulting expression. That the coefficients $B_i$ are as
claimed in the appendix is a tedious but routine computation that
we performed with the help of {\sc Maple}. In particular, the
coefficients of $X$ and $X^3$ vanish identically in $y$ (or in
$t$). The $B_i$ are analytic since they are elementary functions
of the~$D_i$.
\end{proof}

%El fet que $B_1=0$ prove del fet que $D_1$ de BGW es 0. Pero per
%provar que $B_3=0$ treballem amb $B'$. Donat que $B$ es
%essencialment $\int_o^D L ds$, tenim que
% $$
%    B' === \int_0^D L'ds + D' L(x,y,D).
%$$
%Pero per definicio de $D$, $L(x,y,D)=0$. I hem calculat $B'$ fent
%aquesta integral, amb metodes semblants als de $\int_0^D Mds$.

%%%%%%%%%%%%%%%%%%%%%%%%%%%%%%%%%
\section{Asymptotic estimates}\label{sec:est}
%%%%%%%%%%%%%%%%%%%%%%%%%%%%%%%%%

In order to prove Theorem \ref{th:main}, first we need to locate
the dominant singularity $\rho = \gamma^{-1}$ of $G(x)$. Since
$G(x) = \exp (C(x))$, the functions $G(x)$ and $C(x)$ have the
same singularities; hence from now on we concentrate on $C(x)$.

We rewrite the second equation in Lemma \ref{le:BCG} as
\begin{equation}\label{eq:main}
    F(x) = x \exp(B'(F(x))),
\end{equation}
where $F(x)=x C'(x)$. Notice that the singularities of $B'(x)$ and
$F(x)$ are the same, respectively, as those of $B(x)$ and $C(x)$.
 From (\ref{eq:main}) it follows that
\begin{equation}\label{eq:psi}
    \psi(u) = u e^{-B'(u)}
\end{equation}
is the functional inverse of $F(x)$. The dominant singularity of
$\psi$ is the same as that of $B(x)$, which according to
Lemma~\ref{sing-B} is equal to~$R = R(1)$. In order to determine
the dominant singularity $\rho$ of $F(x)$, we have to decide which
of the following possibilities hold; see~Proposition~IV.4
in~\cite{FS} for an explanation.
\begin{enumerate}
  \item There exists $\tau \in (0,R)$ (necessarily unique)
  such that $\psi'(\tau)=0$. Then $\psi$ ceases to be invertible at $\tau$ and
  $\rho = \psi(\tau)$.
  \item We have $\psi'(u) \ne 0$ for all $u \in (0,R)$. Then $\rho =
   \psi(R)$.
\end{enumerate}
The condition $\psi'(\tau)=0$ is equivalent to $B''(\tau) =
1/\tau$. Since $B''(u)$ is increasing (the series $B(u)$ has
positive coefficients) and $1/u$ is decreasing, we are in case (2)
if and only if $B''(R) < 1/R$. Next we show that this is the case.

\begin{claim}\label{cla:R}
Let $R$ be as before the radius of convergence of $B(x)$. Then
$B''(R) < 1/R$.
\end{claim}

\begin{proof}
Lemma~\ref{sing-B} implies that $B''(R) = 2B_4/R^2$ (see
(\ref{B'}) below). Hence the inequality becomes $2B_4 < R$. It
holds because $R \approx 0.0381$ and $B_4 \approx 0.000767$.
\end{proof}
Let us remark that in a related problem, counting series-parallel
graphs, a very similar situation appears but the analogous $\psi$
function \emph{does} have a maximum in its domain of definition
\cite{sp}.

We are now ready for the main result.

\medskip
 \noindent{\emph{Proof of Theorem \ref{th:main}.}}
As we have seen in the previous claim, the dominant singularity of
$F(x)$ is at $\rho = \psi(R)$. In order to obtain the singular
expansion of $F(x)$ at $\rho$, we have to invert the singular
expansion of $\psi(u)$ at $R$.

The expansion of $B'(x)$ follows directly by differentiating the
one in Lemma~\ref{sing-B}:
\begin{equation}\label{B'}
  B'(x) = - {1 \over R} \left( B_2 + 2B_4 X^2 + {5 \over 2} B_5
  X^3 \right) +  \mathcal{O} (X^4).
\end{equation}
Because of $\psi(x) = x \exp(-B'(x))$, by functional composition
we obtain
 $$
  \psi(x) = R e^{B_2/R} \left( 1 + \left({2B_4 \over R} -1\right) X^2
  + {5 B_5 \over 2R} X^3 \right) +  \mathcal{O} (X^4).
$$
Since we are inverting at the singularity, $F(x)$ also has a
singular expansion of square-root type
 $$
 F(x) = F_0 + F_1 X + F_2 X^2 + F_3 X^3 +  \mathcal{O} (X^4),
  $$
with the difference that now $X= \sqrt{1 - x/\rho}$. Given that
$F(x)$ and $\psi(x)$ are functional inverses, the $F_i$ can be
found by indeterminate coefficients, and they turn out to be, in
terms of $R$ and the $B_i$,
\begin{equation}\label{Fi}
 F_0 = R, \quad F_1 = 0, \quad F_2 = {R^2 \over 2B_4-R},
 \quad F_3 = -{5\over2}B_5 (1-2B_4/R)^{-5/2}.
\end{equation}

The singular expansion of $C(x)$ is obtained by integrating $C'(x)
= F(x)/x$, and one gets
\begin{equation}\label{singC}
 C(x) = C_0 + C_2 X^2 + C_4 X^4 + C_5 X^5 +
 \mathcal{O} (X^6).
\end{equation}
The $C_i$, except $C_0$,  are computed easily in terms of the
$F_i$ in equation (\ref{Fi}), and they turn out to be
\begin{equation}
 C_2 = -F_0,  \quad C_4 = -{F_0+F_2 \over 2},
 \quad C_5 = -{2\over5}F_3.
\end{equation}

By singularity analysis (Proposition \ref{prop:transfer}), we
obtain the estimate
 $$
    c_n \sim c \cdot n^{-7/2} \rho^{-n} n!,
 $$
where $c = C_5 / \Gamma(-5/2)$.

However, the coefficient $C_0 = C(\rho)$ is indeterminate after
the integration of $F(x)/x$, and is needed later. To compute it,
we start by integrating by parts
 $$
    C(x) = \int_0^x {F(s) \over s} \, ds = F(x) \log x - \int_0^x
    F'(s) \log s \, ds.
 $$
We change variables $t=F(s)$, so that $s = \psi(t)=t e^{-B'(t)}$,
and the last integral becomes
 $$
 \int_0^{F(x)} \log \psi(t)\, dt =
 \int_0^{F(x)} (\log t - B'(t))  \, dt =
 F(x) \log F(x) - F(x) - B(F(x)).
 $$
Hence
 $$
 C(x) = F(x) \log x  - F(x)\log F(x) + F(x) + B(F(x)).
 $$
Taking into account that $F(\rho) = R$ and $B(R) = B_0$, we get
$$C_0 = C(\rho) = R \log \rho  - R \log R  + R + B_0.$$ A simple
computation shows that, equivalently,
\begin{equation}\label{C(rho)}
C_0 = R + B_0 + B_2.
\end{equation}

The final step is simpler since $G(x) = e^{C(x)}$. We apply the
exponential function to (\ref{singC}) and obtain the singular
expansion
\begin{equation}\label{singG}
 G(x) = e^{C_0}\left(1 + C_2 X^2 + (C_4 + {1 \over 2}C_2^2) X^4 + C_5 X^5\right) +
 \mathcal{O} (X^6),
\end{equation}
where again $X= \sqrt{1 - x/\rho}$. Again by singularity analysis,
we obtain the estimate
 $$g_n \sim g \cdot n^{-7/2} \rho^{-n}n!,$$
where $g = e^{C_0}c$.
 Finally, since $\rho = \psi(R) = R e^{-B'(R)}$ and
$B'(R) = -B_2/R$, we get
 $$
    \rho = R e^{B_2/R},
    \qquad \gamma = \rho^{-1} = {1 \over R} e^{-B_2/R}.
 $$
The constants $c, g$ and $\rho$ can be found using the known value
of $R$ and the expressions for the $B_i$ in the appendix; the
approximate values in the statement have been computed using these
expressions.

Notice that the probability that a random planar graph is
connected is equal to $c_n/g_n \sim c/g = e^{-C_0}$. This result
reappears later in Theorem~\ref{th:components}. \qed

%%%%%%%%%%%%%%%%%%%%%%%%%%%%%%%%%
\section{Gaussian limit laws}\label{sec:laws}
%%%%%%%%%%%%%%%%%%%%%%%%%%%%%%%%%

The proofs in this section
 %of Theorems~\ref{th:edges} and \ref{th:blocs}
are based on bivariate singular expansions and perturbation of
singularities.
To simplify the notation, in this section we denote by
$f'(x,y)$ the derivative of a bivariate function with respect to
$x$.

\medskip
 \noindent{\emph{Proof of Theorem \ref{th:edges}.}}
 We rewrite the second
equation in Lemma~\ref{le:BCG-biv} as
\begin{equation}
    F(x,y) = x \exp(B'(F(x,y),y)),
\end{equation}
where $F(x,y)=x C'(x,y)$. It follows that, for $y$ fixed,
\begin{equation}\label{eq:psi-biv}
    \psi(u,y) = u e^{-B'(u,y)}
\end{equation}
is the functional inverse of $F(x,y)$.

We know from the previous section that $\psi'(u,y)$  does not
vanish for $y=1$ and $u \in (0,R)$, and that $\rho = \psi(R)$ is
the dominant singularity of $F(x)$. Hence by continuity the same
is true for $y$ close to~$1$, and the dominant singularity of
$F(x,y)$ is at
\begin{equation}\label{rho}
 \rho(y) = \psi(R(y),y) = R(y) e^{-B'(R(y),y)}.
\end{equation}
 Given the analytic expressions for the functions involved, the
univariate singular expansion of $\psi(x)$ extends to an expansion
of $\psi(x,y)$ for $y$ fixed. The same is true then for $F(x,y)$
and $C(x,y)$, and we obtain a bivariate expansion
 $$
 C(x,y) = C_0(y) + C_2(y) X^2 + C_4(y) X^4 + C_5(y) X^5 +
 \mathcal{O} (X^6),
 $$
where the $C_i(y)$ are analytic functions, and now $X = \sqrt{1 -
x/\rho(y)}$.

Then Proposition \ref{th:quasi-powers} implies a limit normal law
for the number of edges in random connected planar graphs, with
expectation and variance linear in~$n$. The constants $\kappa$ and
$\lambda$ in the statement of Theorem~\ref{th:edges} are given by
$$
  \kappa = -{\rho'(1) \over \rho(1)}, \qquad
  \lambda = -{\rho''(1) \over \rho(1)}  -{\rho'(1) \over \rho(1)}
    + \left( {\rho'(1) \over \rho(1)}\right)^2,
   $$
   where $\rho'(y) = d \rho(y) / dy$.
Since $G(x,y)$ and $C(x,y)$ have the same dominant singularities
$\rho(y)$, the previous statement also holds  for arbitrary planar
graphs, with the same values of $\kappa$ and $\lambda$.

In order to determine the parameters exactly, we need only an
explicit expression for $\rho(y)$. The expansion (\ref{B'})
extends to an expansion of $B'(x,y)$, whose constant term is
$B'(R(y),y) = -B_2(y)/R(y)$. Hence from (\ref{rho}) it follows
that
 $$
 \rho(y) = R(y) \exp\left( B_2(y)/R(y)\right).
 $$
The appendix contains an explicit expression for $\rho(y)=q(t)$ as
a function of $t$. The necessary derivatives are computed as
$\rho'(y) = q'(t) / Y'(t)$, and the same goes for $\rho''(y)$.
% It is clear that $\kappa$ and $\lambda$ can be expressed explicitly
% in terms of $R(y)$ and $B_2(y)$ and their derivatives at $y=1$.
The approximate values in the statement have been computed in this
way.
 \qed

\medskip
 \noindent{\emph{Proof of Theorem \ref{th:local}.}}
Consider the generating function $C_u(x,y)=C(x,uy)$, where $u$ is
a fixed constant. In this situation the singularity $\rho_u(y)$ of
$C_u$ is given by $\rho(uy)$, and the associated probabilities
$p^u_{n,k}$ of $C_u$ are

\begin{equation}\label{eq:local1}
  p^u_{n,k} = \frac{[y^k][x^n]C_u(x,y)}{[x^n]C_u(x,1)}=
\frac{u^k c_{n,k}}{n! [x^n]C(x,u)}.
\end{equation}
In order to apply Proposition \ref{prop:local} to $C_u$ we need to
know the singularities of $C(x,y)$ when $y$ is away from $1$. The
following claim extends Claim~\ref{cla:R} and shows that the
bivariate singularity expansions given in the proof of
Theorem~\ref{th:edges} hold for every $y$.

\begin{claim}
Let $R(y)$ be the radius of convergence of $B(x,y)$ for $y$ fixed. Then
$B''(R(y),y) < 1/R(y)$.
\end{claim}

\begin{proof}
As in the proof of Claim~\ref{cla:R}, it is enough to show that
$2B_4(y)<R(y)$ for $y\in(0,\infty)$; equivalently, that
$2B_4(t)<\xi(t)$ for $t\in(0,1)$. We bound the logarithm that
appears in  the expression for  $B_4$ (see the appendix) as
$$
\log\left(\frac{1+t}{\sqrt{1+2t}}\right)  \leq
\frac{1+t}{\sqrt{1+2t}}-1.
 $$
Let $\widetilde{B}_4$ be the function obtained by substituting the
logarithm in $B_4$ for the right-hand side in the previous
inequality. Then  it is enough to show that
$$  2\widetilde{B}_4(t) < \xi(t) \qquad \hbox{for $t\in(0,1)$}.
$$
Since both $\widetilde{B}_4$ and $\xi$ are rational functions, the
problem reduces to showing that a certain polynomial (in fact, of
degree 20) with integer coefficients has no root in~$(0,1)$. We
have checked that this is indeed the case using {\sc Maple}.
\end{proof}

Another requirement  is that $\rho(z)$ attains uniquely its
minimum on $|z|=u$ at $z=u$. Suppose it exists $w\neq u$ with
$|w|=u$ such that $|\rho(w)|\leq\rho(u)$. It follows from %equation
(\ref{rho}) that $R(z)$ is equal to $F(\rho(z),z)$, and since
$F(x,y)$ has non-negative coefficients, $|R(w)|=|F(\rho(w),w)|\leq
F(\rho(u),u)=R(u)$. However,  this contradicts the fact that
$R(z)$ attains uniquely its minimum on $|z|=u$ at $z=u$, as shown
in~\cite[Lemma~3]{bender}.

Now Proposition \ref{prop:local} applied to $C_u$ yields
\begin{equation}
\label{eq:local2}
  p^u_{n,\lfloor \mu(u) n\rfloor} \sim \frac{1}{\sqrt{2 \pi n}\, \sigma(u)},
\end{equation}
where $\mu(u)$ are $\sigma(u)$ are given by
$$
  \mu(u)=-\frac{\rho_u'(1)}{\rho_u(1)}=-\frac{u\rho'(u)}{\rho(u)},
$$
$$
  \sigma(u)^2=-\frac{\rho_u''(1)}{\rho_u(1)}-\frac{\rho_u'(1)}{\rho_u(1)}+
     \left(\frac{\rho_u'(1)}{\rho_u(1)}\right)^2
   =-u^2\frac{\rho''(u)}{\rho(u)}-u\frac{\rho'(u)}{\rho(u)}
     +\left(u\frac{\rho'(u)}{\rho(u)}\right)^2.
$$

\smallskip\noindent
Theorem \ref{th:local} follows by combining equations
(\ref{eq:local1}) and (\ref{eq:local2}) for $k=\lfloor \mu(u)
n\rfloor$ and using the asymptotic expression of $[x^n]C(x,y)$ for
$y=u$. The value $\mu$ is constrained to the interval $(1,3)$
since $\lim_{u\rightarrow 0} \mu(u)=1$ and $\lim_{u\rightarrow
\infty} \mu(u)=3$.

%$\sqrt(n) p_{n,\lfloor \mu n\rfloor}$ converges to $(
%\sqrt{2 \pi}\sigma)^{-1}$.
%$$
%
%$ generating function is
%$$
%p^u_n(y)=\sum_k p^u_{n,k} y^k=\frac{[x^n]G_u(x,y)}{[x^n]G_u(x,1)}
%$$
%where the probabilities $p^u_{n,k}$ are related to $p_{n,k}$ by
%the equati

\medskip
 \noindent{\emph{Proof of Theorem \ref{th:appear}.}}
Let us recall Equation (\ref{eq:main})
 $$
    F(x) = x \exp(B'(F(x))),
$$
where $F(x) = x C'(x)$ is the generating function of rooted
connected planar graphs.  In order to mark appearances of $H$, we
have to look at the root $r$ of a rooted connected graph $G$, and
the blocks to which it belongs; recall this is encoded in the term
$\exp(B'(F(x)))$. We are interested in the blocks which are equal
to a single edge $rv$, and within these blocks to the situation
where vertex $v$ is substituted by a copy of $H$. In this case we
mark an appearance of $H$ with the secondary variable~$y$. If we
let $f(x,y)$ be the corresponding generating function, then the
previous discussion translates into the equation
\begin{equation}\label{eq:appear}
 f(x,y) = x \exp\left(B'(f(x,y)) + (y-1) {x^h\over h!}    \right).
\end{equation}
 Notice that from the definition of appearances we do not need
to take into account the automorphisms of $H$; if a copy of $H$ is
substituted for vertex $v$, there is only one way to do it once
the labels are selected, hence the term $x^h/h!$.

In fact, $f(x,y)$ is not the \emph{exact} counting series, since
it does not take into account the possibility that the root $r$
belongs to a copy of $H$ that appears in $G$. This can be
accounted for as follows. The generating function of rooted
connected graphs where the root belongs to an appearance of $H$ is
$f(x,y) x^h/(h-1)!$, since the root can appear in any of the $h$
vertices of $H$. Hence the generating function that counts exactly
\emph{all} appearances is
 $$
g(x,y) = f(x,y) + (y-1) {x^h \over (h-1)!} f(x,y).
$$
Since $f(x,y)$ and $g(x,y)$ have the same dominant singularity for
any fixed $y$ it does not matter which one we choose for
singularity  analysis; hence in the rest of the proof we work with
$f(x,y)$, defined through (\ref{eq:appear}).

Equation (\ref{eq:appear}) can be rewritten  as
 $$
 f(x,y) = \zeta(x,y) \exp\left(B'(f(x,y))    \right),
 $$
where $\zeta(x,y) = x \exp((y-1) x^h/h!)$. Comparing the previous
equation with (\ref{eq:main}), it follows that
 $$
    f(x,y) = F(\zeta(x,y)).
    $$
 Given that $\rho$ is the dominant singularity of $F(x)$,
the dominant singularity of $f(x,y)$ for fixed $y$ is the smallest
value $\tau(y)$  satisfying
\begin{equation}\label{sigma}
 \zeta(\tau(y),y) =
 \tau(y) \exp\left((y-1){\tau(y)^h\over h!}\right) = \rho.
 \end{equation}
Clearly $\tau(1) = \rho$. In order to compute $\tau'(y)$, we
differentiate (\ref{sigma}), set $y=1$, and obtain $\tau'(1) = -
\rho^{h+1}/h!$. To compute $\tau''(1)$ we differentiate again and,
after a simple computation, we get
 $$
  -{\tau'(1) \over \tau(1)} = {\rho^h \over h!} \, , \qquad
  -{\tau''(1) \over \tau(1)}  -{\tau'(1) \over \tau(1)}
    + \left( {\tau'(1) \over \tau(1)}\right)^2 =\rho.
   $$

 From the singular expansion of $F(x)$ at $\rho$, we derive a
corresponding bivariate singular expansion of $f(x,y)$ at
$\tau(y)$, and again a normal limit law follows from Proposition
\ref{th:quasi-powers}. As in the previous proof, a large deviation
estimate also follows, and from this we obtain the bounds in
(\ref{prob-appear}); the details are omitted to avoid repetition.

\medskip
 \noindent{\emph{Proof of Theorem \ref{th:blocs}.}}
The proof is similar to the previous proofs, and so we omit some
details. The generating function $C_1(x,y)$ of connected planar
graphs according to the number of vertices and blocks satisfies
the equation
 $$
 x C_1'(x,y) = x \exp\left(y \, B'\left(x C_1'(x,y)\right)
    \right),
 $$
where $B(x)$ is the univariate generating function of 2-connected
planar graphs.

Let $F_1(x,y) =   x  C_1'(x,y)$. Then, for $y$ fixed,
 $$
    \psi_1(u,y) = u e^{-y B'(u)}
 $$
is the functional inverse of $F_1(x,y)$. The dominant singularity
of $\psi_1(u,y)$ is at $R$, which in this case is independent of
$y$, and the dominant singularity of $F_1(x,y)$ is at
 $$
 \rho_1(y) = \psi_1(R,y) = R e^{-y B'(R)}.
$$
Again we have bivariate singular expansions whose coefficients are
analytic functions of $y$, and the quasi-powers theorem implies
asymptotic normality of the parameter. The asymptotic expressions
for the expected value and variance are obtained as before, but in
this case the computations are particularly easy, since
 $$\rho_1'(y) = -\rho_1(y) B'(R).$$
We know that  $\rho = \psi(R) = R e^{- B'(R)}$, hence
 $$
  \zeta = -{\rho'_1(1) \over \rho_1(1)} = B'(R) = \log(R/\rho)
 \approx 0.03905180273.
  $$
A similar computation gives
 $$
    -{\rho''_1(1) \over \rho_1(1)}  -{\rho'_1(1) \over \rho_1(1)}
    + \left( {\rho'_1(1) \over \rho_1(1)}\right)^2 = B'(R) = \zeta.
\qed
$$

%%%%%%%%%%%%%%%%%%%%%%%%%%%%%%%%%
\section{Poisson limit laws}\label{sec:poisson-laws}
%%%%%%%%%%%%%%%%%%%%%%%%%%%%%%%%%

As opposed to the proofs in the previous section, to prove
Theorems \ref{th:components} and \ref{th:family}, univariate
asymptotics is
 enough.

\medskip
 \noindent{\emph{Proof of Theorem \ref{th:components}.}}
 Let $\nu = C(\rho) = C_0$, the evaluation of $C(x)$ at its dominant
 singularity.
For fixed $k$, the generating function of planar graphs with
exactly $k$ connected components is
 $$
    {1 \over k!} C(x)^k.
 $$
For fixed $k$ we have
 $$
    [x^k] C(x)^k  \sim kC_0^{k-1} [x^n] C(x).
 $$
Hence the probability that a random planar graphs has exactly $k$
components is asymptotically
$$
{ [x^n] C(x) / k! \over [x^n]G(x)} \sim {k C_0^{k-1} \over k!} \,
e^{-C_0} = {\nu^{k-1} \over (k-1)!} \, e^{-\nu},
$$
as was to be proved.
 \qed

\medskip
 \noindent{\emph{Proof of Theorem \ref{th:family}.}}
The proof is similar to the previous one. The generating function
of planar graphs with no component belonging to $\mathcal A$ is
$\exp(C(x)-A(x))$. Hence the generating function of planar graphs
with exactly $k$ components in $\mathcal A$ is
 $$
    {1 \over k!} A(x)^k \exp(C(x)-A(x)) = {1 \over k!} A(x)^k
    e^{-A(x)} G(x).
 $$
The same kind of simple calculation as before gives that the
probability that  a random planar graphs has exactly $k$
components in $\mathcal A$ is asymptotically
 $$
 {A(\rho)^k \over k!} \, e^{-A(\rho)}.
 $$
This finishes the proof of the theorem. \qed

\section{Concluding remarks}

We have found a solution to the problem of counting labelled
planar graphs; however, counting \emph{unlabelled} planar graphs
appears to be much more difficult. If $u_n$ is the number of
unlabelled planar graphs on $n$ vertices, then it is known that
the following limit exists
 $$
    \gamma_u = \lim_{n\to\infty} (u_n)^{1/n},
    $$
and that $\gamma < \gamma_u$, where $\gamma$ is as in
Theorem~\ref{th:main} (see \cite{dominic}). The reason for the
strict inequality $\gamma < \gamma_u$ is that, contrary to what
happens for unrestricted graphs, a planar graph has with high
probability an exponential number of automorphisms~\cite{dominic}.

The best upper bound obtained so far is $\gamma_u < 30.06$. This
is proved in~\cite{upper2} by showing that an unlabelled planar
graph with $n$ vertices can be encoded with $4.91 n$ bits. On the
other hand, our determination of $\gamma $ provides a lower bound
on $\gamma_u$, and shows that at least $4.76 \approx \log_2
\gamma$ bits per vertex are needed.

 We believe that to determine $\gamma_u$ exactly is a very hard
problem, not to speak of determining the subexponential behavior
of $u_n$. The reason is that the equations connecting the
generating functions of labelled planar graphs with different
connectivity requirements, do not hold anymore in the unlabelled
case.

A related problem  is to estimate the number of planar graphs with
a given number of \emph{edges}. In \cite{upper2} it is proved that
an unlabelled planar graph with $m$ edges can be encoded with
$2.82 m$ bits. We can show that at least $2.59$ bits per edge are
needed, as follows.

The coefficient of $y^m$ in $G(1,y)$  is equal to
 $$
 h_m = \sum_n {g_{n,m} \over n!},
 $$
where $g_{n,m}$ is the number of labelled planar graphs with $n$
vertices and $m$ edges. Since a graph on $n$ vertices has at most
$n!$ automorphisms, the number of unlabelled planar graphs with
$m$ edges is at least $h_m$.

The exponential growth of the $h_m$ is determined by the smallest
singularity $\tau$ of $G(1,y)$. Since the smallest singularity of
$G(x,y)$ for fixed $y$ is $\rho(y)$, as given in~(\ref{rho}), it
follows that $\tau$ is the smallest solution to $\rho(\tau) = 1$.
It can be computed exactly with the expressions in the appendix
and it turns out that
 $$
 \lim_{n \to \infty} (h_m)^{1/m} =\tau \approx 6.03 \approx 2^{2.59}.
 $$

Finally, let us mention that the explicit expressions we have
obtained for the generating functions of labelled planar graphs
have been applied to the design of very efficient algorithms for
generating random planar graphs uniformly \cite{eric}.

%\newpage
%%%%%%%%%%%%%%%%%%%%%%%%%%%%%%%%%
\section*{Appendix}\label{sec:app}
%%%%%%%%%%%%%%%%%%%%%%%%%%%%%%%%%

Here we list the functions $B_0, B_2, B_4, B_5$ that have been
used in the previous sections, as functions of $t$. As has been
explained already, to become functions of $y$ near $1$, they must
be evaluated at the unique solution of $Y(t) = y$. For
completeness, we also write down the function $\rho$.

\begin{eqnarray*}
B_0 &=& {\frac { \left( 3\,t-1 \right) ^{2} \left( 1+t
 \right) ^{6}\log  \left( 1+t \right) }{512\,{t}^{6}}}\,-\,{\frac {
 \left( 3\,{t}^{4}-16\,{t}^{3}+6\,{t}^{2}-1 \right) \log  \left( 1+3\,t
 \right) }{32\,{t}^{3}}} \\
 && -\,{\frac { \left( 1+3\,t
 \right) ^{2} \left( 1-t \right) ^{6}\log  \left( 1+2\,t \right) }{1024\,{t}
^{6}}}+ {1\over 4}\,\log  \left( 3+t \right) -{1\over2}\,\log
\left( t \right)
-{3\over8}\, \log  \left( 16 \right)\\
 && -\,{\frac {
\left( 217\,{t}^{6
}+920\,{t}^{5}+972\,{t}^{4}+1436\,{t}^{3}+205\,{t}^{2}-172\,t+6
 \right)  \left( 1-t \right) ^{2}}{2048\,{t}^{4} \left( 1+3\,t \right)
 \left( 3+t \right) }}
\end{eqnarray*}

\bigskip

\begin{eqnarray*}
B_2 &=& {\frac { \left( 1-t \right) ^{3} \left( 3\,t-1
 \right)  \left( 1+3\,t \right)  \left( 1+t \right) ^{3}\log  \left( 1+
t \right) }{256 \, {t}^{6}}}\\
 &&-{\frac { \left( 1-t \right) ^{3}
\left( 1 +3\,t \right) \log  \left( 1+3\,t \right)
}{32\,{t}^{3}}}+{ } \,{\frac { \left( 1+3\,t \right) ^{2} \left(
1-t \right) ^{6}\log
 \left( 1+2\,t \right) }{512\,{t}^{6}}}\\
 &&+\,{\frac { \left(
1-t \right) ^{4} \left(
185\,{t}^{4}+698\,{t}^{3}-217\,{t}^{2}-160\,t +6 \right)
}{1024\,{t}^{4} \left( 1+3\,t \right)  \left( 3+t \right) }}
\end{eqnarray*}

\bigskip

%\begin{eqnarray*}
%B_4 &=& \,{\frac {\log  \left( 1+t \right)  \left( -1+t
% \right) ^{6} \left( 1+3\,t \right) ^{2}}{512\,{t}^{6}}}-
%\,{\frac { \left( 1+3\,t \right) ^{2} \left( -1+t \right) ^{6}\log
% \left( 1+2\,t \right) }{{1024\,t}^{6}}}\\
% &&-\,{\frac { P \, \left( -1+t \right)
%^{5}}{2048\,{t }^{4} \left( 1+3\,t \right) ^{2} \left( 3+t \right)
%^{2} Q } }
%\end{eqnarray*}

\begin{eqnarray*}
B_4 &=& \,{\frac {\log  \left({\displaystyle \frac{1+t}{\sqrt{1+2\,t}}} \right)  \left( 1-t
 \right) ^{6} \left( 1+3\,t \right) ^{2}}{512\,{t}^{6}}}
+\,{\frac { P \, \left( 1-t \right)
^{5}}{2048\,{t }^{4} \left( 3+t \right) Q } }
\end{eqnarray*}

\bigskip

\begin{eqnarray*}
%B_5 &=&-{\frac {\sqrt {3}}{90}} ( -1+t ) ^6 \left( {\frac {
% ( 1+3t )  ( 3+t ) {T}} { ( 1+t) t^2 Q}} \right) ^{3/2} ( 1+3t
% )^{-4} \qquad \qquad \qquad \\
% &&
%+ {\frac {4}{45}} {\frac { ( 18+6t+t^2 ) (-1+t)^7 (1+t)^{2}
%(1+2t)^2 S^{3/2} \sqrt{3}} { (1+3t) (t (1+t) Q ) ^{5/2} } }
B_5 & = & -\frac{\sqrt{3}}{90} \frac{(1-t)^6}{(1+t)^{3/2}} \left( \frac{S}{t\,Q}\right)^{5/2}
\end{eqnarray*}

\bigskip \noindent  where

\begin{eqnarray*}
P &=&-2400 +
57952\,{t} +303862\,{t}^2 +466546\,{t}^3 +
264775\,{t}^4 +76679\,{t}^5 +11495\,{t}^6+
739\,{t}^7\\
%5627\,{t}^{10}+117602\,{t}^{9}+1034079\,{t}^{8}+5055924\,{t}^{7}+
%14504931\,{t}^{6}+23910962\,{t}^{5}\\
%&& +19869793\,{t}^{4}+8722840\,{t}^{3}
%+2071714\,{t}^{2}+211200\,t-7200 \\
 Q &=&
400+1808\,t+2527\,{t}^{2}+1155\,{t}^{3}+237\,{t}^{4}+17\,{t}^ {5}
\\
S & = & 144+592\,t+664\,t^2+135\,t^3+6\,t^4-5\,{t}^5
%S &=& 144+592\,t+664\,{t}^{2}+135\,{t}^{3}+6\,{t}^{4}-5\,{t}^{5}\\
%T &=&
%17\,{t}^{6}+154\,{t}^{5}+469\,{t}^{4}+376\,{t}^{3}+280\,{t}^{
%2}+192\,t+48
\end{eqnarray*}

%\newpage
\bigskip
\noindent $$ \rho \quad = \quad  {1\over 16}\,\sqrt {1+3\,t}
\left( 1-t \right) ^{3}{t}^{-3} \exp(A),$$

\medskip where

\begin{eqnarray*}
A &=&  \, \frac {\log
 ( 1+t)( 3t-1) ( 1+t) ^3}
 { 16\, t^3}   +\,\frac {\log  ( 1+2t ) ( 1+3t) (1-t) ^3}
 {32\,t^3} \\
 && +  \frac {( 1-t )  (185t^4+698t^3-217t^2- 160t+6) }{64\,t
( 1+3t) ^2 ( 3+t ) } \end{eqnarray*}

\bigskip
\noindent The approximate values of the univariate constants $B_i
= B_i(t)$, where $t=0.6263716633$ is the unique solution of
$Y(t)=1$, are
\begin{eqnarray*}
 && B_0 = 0.7396995711 \cdot 10^{-3}, \quad
 B_2 =-0.1491431215 \cdot 10^{-2}, \\
 &&  B_4 =0.7671782851 \cdot 10^{-3}, \quad
 B_ 5 =-0.3501857790  \cdot 10^{-5}.
\end{eqnarray*}

Finally we include the singular expansion of $U(x,D(x,y))$ at the
dominant singularity $R(y)$:
 $$
 U(x,D(x,y)) = U_0(y) + U_1(y) X + U_2(y) X^2 + \mathcal{O}(X^3),
 $$
where $X = \sqrt{1-x/R(y)}$ and the $U_i$ are, as functions of
$t$, given by
\begin{eqnarray*}
 U_0 &=& \frac{1}{3t} \\
 U_1 &=&
\left(
\frac{4(1+3t)^2(-5t^5+6t^4+135t^3+664t^2+592t+144)}{27t^3(1+t)Q}
\right)^{1/2}
 \\
  U_2 &=& \frac{2(1+3t)T}
 {27t^2(1+t)^2 Q^2}
\end{eqnarray*}
where $Q$ is as before and
\begin{eqnarray*}
T &=&
691t^{12}+10112t^{11}+98693t^{10}+719346t^9+3723625t^8+13180580t^7
  +31133003t^6 \\
&&+
47691938t^5+47354348t^4+30156200t^3+11835336t^2+2596736t+243072
\end{eqnarray*}

\bibliographystyle{amsplain}

\end{document}